\documentclass[11pt]{amsart}
\usepackage{enumerate,amsmath,amssymb,amsfonts,bm}
\usepackage{frcursive,fullpage,graphicx,psfrag}

\newtheorem{Theorem}{Theorem}
\newtheorem{Lemma}{Lemma}

\begin{document}

\author{Abdelmalek Abdesselam}
\address{Abdelmalek Abdesselam,
Department of Mathematics,
P. O. Box 400137,
University of Virginia,
Charlottesville, VA 22904-4137, USA}
\email{malek@virginia.edu}

\title{The weakly dependent strong law of large numbers revisited}

\begin{abstract}
We give a short, self-contained, and elementary proof of the strong law of large numbers under a power law decay hypothesis for joint second moments. The result is related to the classical one by Lyons. However, we also provide a rate of convergence. Our proof does not use maximal inequalities and is instead inspired by the method of multiscale large versus small field decompositions in constructive quantum field theory.
\end{abstract}

\maketitle

\section{Introduction and main theorem}

Let $(X_n)_{n\in\mathbb{N}}$ be a sequence of centered, real-valued, square integrable random variables on the same probability space $(\Omega,\mathcal{F},\mathbb{P})$. We will denote the average of the first $N$ variables by $A_N=\frac{X_1+\cdots+X_N}{N}$.
The main result of this article is as follows.

\begin{Theorem}
Suppose the sequence satisfies
\[
\exists \gamma>0, \exists K>0, \forall (m,n)\in\mathbb{N}^2,\ \ \mathbb{E}X_mX_n\le\frac{K}{(1+|m-n|)^{\gamma}}\ .
\]
Let $\beta$ be a parameter in the interval $\left(\frac{1}{2},1\right)$ if $\gamma\ge 1$,
or in the interval $\left(1-\frac{\gamma}{2},1\right)$ if $0<\gamma<1$. Then, with probability one, we have
\[
|A_N|=O\left(\frac{\log n}{n^{1-\beta}}\right)\ .
\] 
\end{Theorem}
Note that our hypothesis (for $m\neq n$) automatically holds in the case of negatively correlated variables.
Also note that our hypothesis includes (when $m=n$) the requirement of uniformly bounded variances for the $X_n$, just as in~\cite[Corollary 11]{Lyons}. The result by Lyons allows more general bounds (there denoted by $\Phi_1$) on second moments. However, in the power law case (i.e., $\Phi_1(x)=(1+x)^{-\gamma}$) the corresponding hypothesis, namely $\gamma>0$, is identical to ours.
On the other hand,~\cite{Lyons} proves the strong law of large numbers (SLLN), $A_N\rightarrow 0$ a. s., yet without rate of convergence such as the one provided by our theorem. The SLLN for dependent random variables
(with or without almost sure rate of convergence) has been investigated in a number of relatively recent articles. In addition to~\cite{Lyons}, see for
instance~\cite{ChandraG,FazekasK,Kuczmaszewska,YangSY,HuRV,Sung,DoukhanKL,Korchevsky1,Korchevsky2}.
For example,~\cite[Theorem 2]{Korchevsky2} implies the $\gamma>1$ part of our theorem but does not cover the case of long-range dependence $0<\gamma<1$. In this note, we did not aim for maximal generality but rather for maximal simplicity. Indeed, in most of the literature we cited, the SLLN is proved by a two-step procedure where the intermediate stage consists in establishing a suitable maximal inequality. Our proof, inspired by the multiscale large versus small field decomposition method in constructive quantum field theory (see, e.g.,~\cite{AbdesselamR}), is direct and bypasses the need for maximal inequalities. 
It is based on two simple ingredients. The first one is what one may loosely call multiscale (or dyadic) analysis, i.e., studying a random function (here $n\mapsto X_n$) in terms of its sums or averages on dyadic blocks. The latter are most easily visualized thanks to a dyadic tree. The second ingredient is combinatorial optimization in order to get good estimates. This involves the use of a very simple algorithm, namely, the greedy algorithm which can be summarized by the phrase ``grab as much as you can, as soon as you can''.

\section{Proof of the theorem}

For $\gamma>0$, let us define the nondecreasing function $p_\gamma:[1,\infty)\rightarrow(0,\infty)$ as follows.
\[
p_{\gamma}(x)=\left\{
\begin{array}{cll}
\frac{1}{\gamma-1} & {\rm if} & \gamma>1\ , \\
(\log x)+1 & {\rm if} & \gamma=1\ , \\
\frac{x^{1-\gamma}}{1-\gamma} & {\rm if} & \gamma<1\ . 
\end{array}
\right.
\]

\begin{Lemma}\label{sumLem}
For every finite set $J\subset\mathbb{N}$,
\[
\sum_{n\in J}\frac{1}{(1+n)^{\gamma}}\le p_{\gamma}(|J|+1)
\]
where $|J|$ denotes the cardinality of $J$.
\end{Lemma}

\noindent{\bf Proof:}
The inequality is trivial if $J=\varnothing$. Otherwise the left-hand side is maximized, for fixed $|J|$, when $J=\{1,\ldots,|J|\}$. A simple sum/integral comparison gives the upper bound $\int_0^{|J|}\frac{dx}{(1+x)^{\gamma}}$ and the lemma follows from the evaluation of the integral in all three cases for $\gamma$. \qed

\begin{Lemma}\label{sumsumLem}
For every nonempty set $J\subset\mathbb{N}$,
\[
\sum_{(m,n)\in J^2}\frac{1}{(1+|m-n|)^{\gamma}}\le |J|\times(1+2p_{\gamma}(|J|))\ .
\]
\end{Lemma}

\noindent{\bf Proof:}
By separating the cases corresponding to the relative positions of $m$ and $n$, one sees that the left-hand side is equal to
\[
|J|+2\sum_{m\in J}\sum_{k\in J_m}\frac{1}{(1+k)^{\gamma}}
\]
where $J_m=\{n-m\ |\ n\in J\ {\rm and}\ n>m\}$. We then apply Lemma \ref{sumLem} to the sum over $J_m$, together with the nondecreasing property of $p_{\gamma}$ and the obvious inequality $|J_m|+1\le |J|$, in order to conclude.
\qed

For any nonempty finite set $J\in \mathbb{N}$, we have, by hypothesis,
\[
\mathbb{E}\left(\sum_{n\in J}X_n\right)^2=\sum_{(m,n)\in J^2}\mathbb{E}X_m X_n\le K
\sum_{(m,n)\in J^2}\frac{1}{(1+|m-n|)^{\gamma}}\ .
\]
For every $c>0$, we have, using Chebychev's Inequality and Lemma \ref{sumsumLem},
\begin{equation}
\mathbb{P}\left(\left|\sum_{n\in J}X_n\right|\ge c\right)\le K c^{-2}|J|\times(1+2p_{\gamma}(|J|))\ .
\label{ChebychevEq}
\end{equation}
For every $(k,i)\in\mathbb{N}_0\times\mathbb{N}$, we define the dyadic block 
\[
B_{k,i}=\left\{
n\in\mathbb{N}\ |\ (i-1)2^k+1\le n\le i 2^k
\right\}\ .
\]
It is convenient to visualize them using an infinite tree as in the figure:
\[
\parbox{12cm}{
\psfrag{1}{$\scriptstyle{1}$}
\psfrag{2}{$\scriptstyle{2}$}
\psfrag{3}{$\scriptstyle{3}$}
\psfrag{4}{$\scriptstyle{4}$}
\psfrag{5}{$\scriptstyle{5}$}
\psfrag{6}{$\scriptstyle{6}$}
\psfrag{7}{$\scriptstyle{7}$}
\psfrag{8}{$\scriptstyle{8}$}
\psfrag{9}{$\scriptstyle{9}$}
\psfrag{a}{$\scriptstyle{10}$}
\psfrag{b}{$\scriptstyle{11}$}
\psfrag{c}{$\scriptstyle{12}$}
\psfrag{d}{$\scriptstyle{13}$}
\psfrag{e}{$\scriptstyle{14}$}
\psfrag{f}{$\scriptstyle{15}$}
\psfrag{g}{$\scriptstyle{16}$}
\psfrag{p}{$k=0$}
\psfrag{q}{$k=1$}
\psfrag{r}{$k=2$}
\psfrag{s}{$k=3$}
\psfrag{t}{$k=4$}
\psfrag{u}{$\cdots$}
\psfrag{v}{$\ddots$}
\psfrag{w}{$\rightarrow\ i$}
\raisebox{1ex}{
\includegraphics[width=12cm]{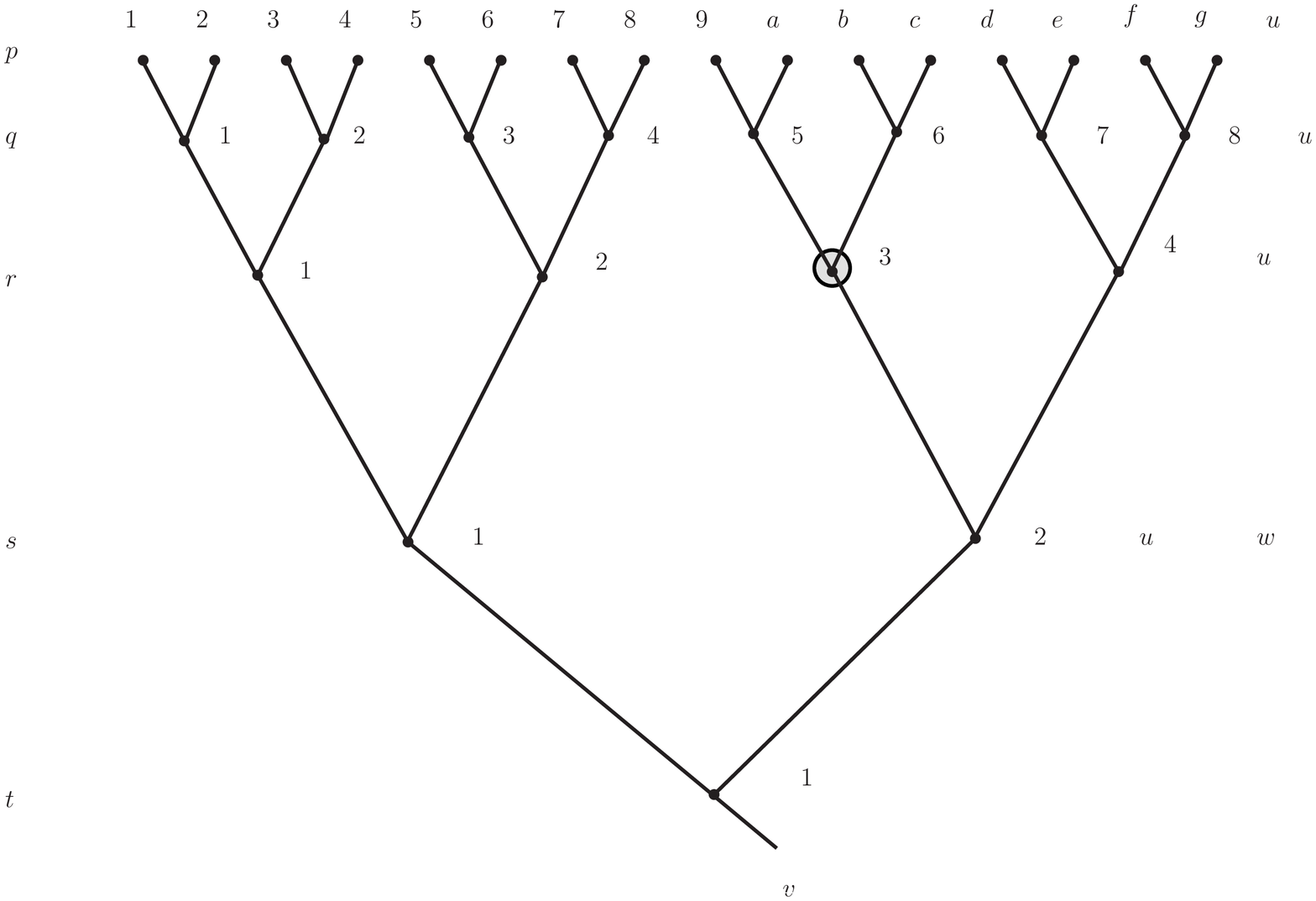}}
}
\]
For example, the node circled in grey with coordinates $(k,i)=(2,3)$ corresponds to the block $B_{2,3}=\{9,10,11,12\}\subset\mathbb{N}$. The numbers indicated on the tree refer to the horizontal $i$ coordinate. The $k$ coordinate indicates the depth. Finally, the set of leafs of the tree is a visualization of the set $\mathbb{N}$ which labels the random variables $X_n$. 
Depending on the realized sample $\omega\in\Omega$, we will call $B_{k,i}$ a bad block (or large field block) if
\[
\left|\sum_{n\in J}X_n\right|\ge 2^{\beta k} i^{\beta}\ .
\] 
Otherwise, we say $B_{k,i}$ is a good block (or small field block).
By (\ref{ChebychevEq}) with $c=2^{\beta k} i^{\beta}$, we have $\forall(k,i)\in\mathbb{N}_0\times\mathbb{N}$,
\[
\mathbb{P}\left(B_{k,i}\ {\rm is\ bad}\right)\le\ K\ i^{-2\beta}\times 2^{(1-2\beta)k}(1+2p_{\gamma}(2^k))
\]
and therefore
\[
\sum_{(k,i)\in\mathbb{N}_0\times\mathbb{N}}\mathbb{P}\left(B_{k,i}\ {\rm is\ bad}\right)<\infty\ .
\]
Indeed, the sum over $i$ converges by the hypothesis $\beta>\frac{1}{2}$. The sum over $k$ is also convergent as can easily be checked in all three cases for $\gamma$. For instance, in the long-range dependence case when $0<\gamma<1$, bounding $\sum_{k\ge 0} 2^{(1-2\beta)k}(1+2p_{\gamma}(2^k))$ amounts to bounding $\sum_{k\ge 0} 2^{(1-2\beta)k}\ 2^{(1-\gamma)k}<\infty$, because of the assumtion $\beta>1-\frac{\gamma}{2}$.
By the first Borel-Cantelli Lemma, it is thus immediate that the (random) set $F\subset \mathbb{N}_0\times\mathbb{N}$ of bad block labels is almost surely finite.

Assuming finiteness of $F$, let $N_F=\max\left(\cup_{(k,i)\in F} B_{k,i}\right)\in\{-\infty\}\cup\mathbb{N}$. The theorem is then a consequence of the following observation.

\begin{Lemma}
If $N\ge 4N_F-1$, then
\[
\left|\sum_{n=1}^{N}X_n\right|\le \left(\frac{\log N}{\log 2}+1\right)\times N^{\beta}\ .
\]
\end{Lemma}

\noindent{\bf Proof:}
Note that $N$ can be uniquely written as $N=2^{k_1}+\cdots+2^{k_l}$ with $k_1>\cdots >k_l\ge 0$. The $k$'s correspond to the positions of the ones in the binary representation of $N$.
Define
\[
\begin{array}{lll}
i_1 & = & 1 \\
i_2 & = & 2^{k_1-k_2}+1 \\
i_3 & = & 2^{k_1-k_3}+2^{k_2-k_3}+1 \\
 & \vdots & \\
i_l & = & 2^{k_1-k_l}+2^{k_2-k_l}+\cdots+ 2^{k_{l-1}-k_l}+1\ .
\end{array}
\]
Then $B_{k_1,i_1},\ldots, B_{k_l,i_l}$ form a set partition of $\{1,2,\ldots, N\}$. It is the partition provided by the greedy algorithm, namely, $B_{k_1,i_1}=\{1,2,\ldots, 2^{k_1}\}$ is the biggest dyadic block inside $\{1,2,\ldots, N\}$ and starting from $1$, while $B_{k_2,i_2}=\{2^{k_1}+1,\ldots,2^{k_1}+2^{k_2}\}$ is next biggest one can form, etc.

Provided all the blocks $B_{k_1,i_1},\ldots, B_{k_l,i_l}$ are good, one can write the estimates
\[
\left|\sum_{n=1}^{N}X_n\right|=\left|
\sum_{s=1}^l\sum_{n\in B_{k_s,i_s}}X_n\right|
\le \sum_{s=1}^l\left|\sum_{n\in B_{k_s,i_s}}X_n\right|
\]
\[
\le  \sum_{s=1}^l 2^{\beta k_s} i_s^{\beta}=\sum_{s=1}^{l}
\left(2^{k_1}+2^{k_2}+\cdots+ 2^{k_s}\right)^{\beta}\le l N^{\beta}
\]
by construction.
Since $k_1>\cdots>k_l\ge 0$, we have $l\le k_1+1$. But $2^{k_1}\le N$, so we obtain $l\le \frac{\log N}{\log 2}+1$
and the desired inequality follows.

All that remains is to show that the hypothesis $N\ge 4N_F-1$ is enough to guarantee that all the blocks $B_{k_1,i_1},\ldots, B_{k_l,i_l}$ are good.
This is essentially a geometric argument based on the dyadic tree. 
If $F=\varnothing$, then $N_F=\infty$ and the hypothesis $N\ge 4N_F-1$ is moot. However, in that case there is nothing more to prove since all blocks are good. We now assume $F\neq\varnothing$ (and of course finite).
The condition $N\ge 4N_F-1$ and the greedy algorithm chosen for the construction of $B_{k_1,i_1},\ldots, B_{k_l,i_l}$ ensure that all the bad blocks are strict subsets of $B_{k_1,i_1}$. Indeed, let $2^r$ be the smallest power of two such that $2^r\ge N_F$. Thus all bad blocks should be subsets of $B_{r,1}$. By construction $2^{k_1}\le N<2^{k_1+1}$ and therefore $4N_F\le N+1\le 2^{k_1+1}$. From $N_F\le 2^{k_1-1}$, we deduce $2^r\le 2^{k_1-1}$, i.e., $r+1\le k_1$. Since all bad blocks are strict subsets of $B_{k_1,i_1}$,
none of the blocks $B_{k_1,i_1},\ldots, B_{k_l,i_l}$ can be bad and we are done.

\qed


\begin{thebibliography}{999}

\bibitem{AbdesselamR}
A. Abdesselam and V. Rivasseau,
An explicit large versus small field multiscale cluster expansion.
Rev. Math. Phys. {\bf 9} (1997), no. 2, 123--199. 

\bibitem{ChandraG}
T. K. Chandra and S. Ghosal,
Extensions of the strong law of large numbers of Marcinkiewicz and Zygmund for dependent variables.
Acta Math. Hungar. {\bf 71} (1996), no. 4, 327-–336. 

\bibitem{DoukhanKL}
P. Doukhan, O. Klesov and G. Lang,
Rates of convergence in some SLLN under weak dependence conditions.
Acta Sci. Math. (Szeged) {\bf 76} (2010), no. 3--4, 683-–695. 

\bibitem{FazekasK}
I. Fazekas and O. Klesov,
A general approach to the strong laws of large numbers.
Theory Probab. Appl. {\bf 45} (2002), no. 3, 436--449

\bibitem{HuRV}
T.-C. Hu, A. Rosalsky and A. Volodin,
On convergence properties of sums of dependent random variables under second moment and covariance restrictions.
Statist. Probab. Lett. {\bf 78} (2008), no. 14, 1999--2005. 

\bibitem{Korchevsky1}
V. M. Korchevsky,
On the strong law of large numbers for sequences of random variables without the independence condition.
Vestnik St. Petersburg Univ. Math. {\bf 44} (2011), no. 4, 268--271. 

\bibitem{Korchevsky2}
V. M. Korchevsky,
On the strong law of large numbers for sequences of dependent random variables with finite second moments.
J. Math. Sci. (N.Y.) {\bf 206} (2015), no. 2, 197--206.

\bibitem{Kuczmaszewska}
A. Kuczmaszewska,
The strong law of large numbers for dependent random variables.
Statist. Probab. Lett. {\bf 73} (2005), no. 3, 305--314. 

\bibitem{Lyons}
R. Lyons,
Strong laws of large numbers for weakly correlated random variables.
Michigan Math. J. {\bf 35} (1988), no. 3, 353--359. 

\bibitem{Sung}
S. H. Sung,
Maximal inequalities for dependent random variables and applications.
J. Inequal. Appl. (2008), Art. ID 598319, 10 pp. 

\bibitem{YangSY}
S. Yang, C. Su and K. Yu, Keming,
A general method to the strong law of large numbers and its applications.
Statist. Probab. Lett. {\bf 78} (2008), no. 6, 794--803. 

\end{thebibliography}
\end{document}